\documentclass[12pt]{article}

\usepackage{amsmath,amssymb,amsfonts,theorem,makeidx,latexsym,epsfig,,subfigure}

\newtheorem{defn}{Definition}[section]

\newtheorem{lemma}[defn]{Lemma}

\newtheorem{proposition}[defn]{Proposition}

{\theorembodyfont{\rmfamily}

\newtheorem{ex}[defn]{Example}}

\newtheorem{thm}[defn]{Theorem}

\newtheorem{prop}[defn]{Proposition}

\newtheorem{cor}[defn]{Corollary}

\numberwithin{equation}{section}

\newcommand{\h}{{\mathcal H}}

\newcommand{\mn}{\mathbb N}

\def\bp{{\noindent\bf Proof. \ }}

\def\ep{\hfill$\square$\par\bigskip}

\def\bqs{\begin{equation}}

\def\eqs{\tag*{$\square$}\end{equation}\par\bigskip}

\def\la{\langle}

\def\ra{\rangle}

\def\bop{\begin{op}\rm}

\def\eop{\end{op}}

\def\bee{\begin{eqnarray}}

\def\ene{\end{eqnarray}}

\def\bes{\begin{eqnarray*}}

\def\ens{\end{eqnarray*}}

\def\bei{\begin{itemize}}

\def\eni{\end{itemize}}

\def\bt{\begin{thm}}

\def\et{\end{thm}}

\def\bc{\begin{cor}}

\def\ec{\end{cor}}

\def\bpr{\begin{prop}}

\def\epr{\end{prop}}

\def\bl{\begin{lemma}}

\def\el{\end{lemma}}

\def\bd{\begin{defn}}

\def\ed{\end{defn}}

\def\bex{\begin{ex}}

\def\enx{\end{ex}}

\def\bfi{\begin{fig}}

\def\efi{\end{fig}}

\def\suj{\sum_{j\in J}}

\def\A{\cal A}

\title{Modular Riesz bases versus Riesz bases in Hilbert $C^*$-Modules}

\date{\today}

\author{Marzieh Hasannasab}

\begin{document}

\maketitle

\begin{abstract}
	
	In this paper we give new characterizations of modular Riesz bases in Hilbert $C^*$-modules. We prove that modular Riesz bases share many properties  with Riesz bases in Hilbert spaces. Moreover we show that there are also important differences; for example, there exist  exact frames  that are not modular Riesz bases.
\end{abstract}

\section{Introduction}



The theory of frames in Hilbert $C^*$-modules was introduced by Frank and Larson \cite{frank};
more recent works include   \cite{bakic, larson2,  han-jing, jing,  ostad,  luef}. Hilbert $C^*$-modules  are generalizations of Hilbert spaces in which the inner product takes values in a $C^*$-algebra.
Kasparov's Stabilization Theorem \cite{kas} shows  that
every finitely or countably generated Hilbert $C^*$-module has a Parseval frame.  The differences between a Hilbert space and a Hilbert $C^*$-module have a significant impact  on frame theory in these spaces too. In this paper, we focus on  Riesz bases in Hilbert $C^*$-modules. In the literature, two attempts
to define Riesz bases in Hilbert 
$C^*$-modules have been given.
In \cite{frank}, a Riesz basis in a Hilbert $C^*$-module is defined as an $\omega$-independent frame; in \cite{ostad}, modular
Riesz bases are defined by a direct generalization of one of the equivalent definitions of Riesz bases in Hilbert spaces (see the end of the section for the exact definition). 
We show that modular Riesz bases in Hilbert $C^*$-modules and Riesz bases in Hilbert spaces share many key features; indeed, 
we prove that some of the well-known characterization theorems for Riesz bases in Hilbert spaces have an analog in Hilbert $C^*$-modules.  Furthermore, characterizations of modular Riesz bases with respect to the canonical dual frame and $\omega$-independent frames are  given. As the last result we show that  even though modular Riesz bases behave similar to Riesz bases in Hilbert spaces, there are some important differences due to the special structure of the Hilbert $C^*$-modules. Indeed we show that  Riesz bases  and in particular modular Riesz bases in Hilbert $C^*$-modules are exact frames but there exist exact frames in Hilbert $C^*$-modules that are not modular Riesz bases.


The paper is organized as follows. In the rest of the introduction
we will recall some basic definitions concerning Hilbert $C^*$-modules,  their frames and the two definitions of Riesz bases in such spaces.  The new results appear in Section 2. 

 Throughout  this paper, $\A$ denotes a  unital $C^*$-algebra  with identity $1_{\A}$. An element $a\in\A$ is positive if $a=b^*b$ for some $b\in \A$. Considering a vector space $\h$, assume that there exists an ${\cal A}$-valued inner product $\langle \cdot,\cdot \rangle : \h\times \h\rightarrow \A$ with the following properties:
\begin{enumerate}
\item $\langle x,x\rangle \geq 0$ for every $x\in \h$ and $\langle x,x\rangle=0$ if and only if $x=0$,
\item $\langle x,y\rangle=\langle y,x\rangle^*$ for every $x,y\in \h$,
\item $\langle ax+y ,z\rangle=a\langle x,z\rangle+\langle y,z\rangle$ for every $a\in \A$ and $x,y,z\in \h$.
\end{enumerate}
Under the conditions stated above, $\h$  is called a {\it Hilbert $\A$-module} (or {\it Hilbert $C^*$-module} in case we want to
avoid the reference to the name of the underlying $C^*$-algebra) if it is complete with respect to the norm
$\|x\|:=\|\langle x,x\rangle\|^{\frac{1}{2}}$. 
A Hilbert {$\A$}-module $\h$ is {\it countably generated} if there exists a countable subset $\{x_j:j\in{\Bbb N}\}$ of $\h$ such that the set of all its finite ${\cal A}$-linear combinations is dense in $\h$; in that case the set  $\{x_j:j\in{\Bbb N}\}$  is called a {\it set of generators}. 

 Let $\h$ and ${\cal K}$ be two Hilbert ${\cal A}$-modules over a $C^*$-algebra ${\cal A}$. A linear operator $T:\h\to{\cal K}$ is called {\it adjointable} if the exists a linear operator $T^*:{\cal K}\to\h$ such that 
 \bes \la Tx, y \ra = \la x ,  T^*y \ra, \quad \text{ for all } x\in \h, y\in{\cal K}.
 \ens
Given any unital $C^*$-algebra, every adjoinable operator is bounded and ${\cal A}$-linear, \cite{lance}.
For every countable index set $ J,$  the {\it standard Hilbert ${\cal A}$-module} is defined by
$$\ell^2( J, {\cal A}):={\left\{ \{a_j\}_{j\in J}\subset {\cal A}: ~\sum_{j\in J}a_j^*a_j ~is~ norm ~convergent~ in~ {\cal A}\right\}}.$$
For each $j\in  J$, letting $\delta_{ij}=0$ if $i\neq j$ and $\delta_{jj}=1$, define $e_j=\{\delta_{ij}1_{\cal A}\}_{i\in J}$. The sequence $\{e_j\}_{j\in J}$ is called the {\it canonical orthonormal basis} for $\ell^2( J, {\cal A})$. 

Following \cite{frank}, we will now give the key definition of frames in Hilbert $C^*$-modules.
\bd
A sequence $\{x_j\}_{j\in J}$  of elements in a Hilbert $C^*$-module $\h$ is said to be a  {\it frame}  if the infinite series $\sum_{j\in J}\langle x,x_j\rangle\langle x_j,x\rangle $ converges in norm for all $x\in\h$ and  there exist two constants
	$0<C\leq D<\infty$ such that
	\begin{equation}\label{frame}
	C\langle x,x\rangle\leq\sum_{j\in J}\langle x,x_j\rangle\langle x_j,x\rangle \leq D\langle x,x\rangle,\quad x\in\h.
	\end{equation}
	Appropriate choices for the numbers $C$ and $D$ are called {\it frame bounds}.
 $\{x_j\}_{j\in J}$ is a {\it Bessel sequence} with bound $D$ if the right-hand side inequality in \eqref{frame} holds. A frame is called {\it tight frame} if we can choose $C=D$ and it is called {\it Parseval frame} if it is a tight frame with bound 1. \ed
 Kasparov's Stabilization Theorem \cite{kas} shows that every  finitely or countably generated Hilbert $C^*$-module over a unital $C^*$-algebra has a frame.
	In inequality \eqref{frame} we are comparing the positive elements in  $\A$.  Aramba\v{s}i\'{c} \cite{arambasic} and Jing in \cite{jing}, independently showed that one can replace \eqref{frame} with two inequalities in terms of the norm of elements. Indeed, it is proved that a sequence $\{x_j\}_{j\in J}\subset\h$ is a frame
if and only if there exist positive constants $C$ and $D$ such that
\[ C\|x\|^2\leq\|\sum_{j\in J}\langle x,x_j\rangle\langle x_j,x\rangle\|\leq D\|x\|^2,\quad x\in \h.\]
For every Bessel sequence $\{x_j\}_{j\in J}$, the operator $T:\h\rightarrow l^2( J, {\cal A})$ defined as 
\bes Tx=\{\langle x,x_j\rangle\}_{j\in J}, \quad x\in \h\ens
is called the {\it analysis operator}.  The operator $T$ is adjointable   and its adjoint, the {\it synthesis operator}, is given by $U\{a_j\}_{j\in J}=\sum_{j\in J}a_j x_j.$ The {\it frame operator } $S:\h\rightarrow \h$, defined as $Sx=U^*U(x)=\sum_{j\in J}\langle x,x_j\rangle x_j$,
is bounded, positive and invertible, see \cite{frank}.
Thus the following {\it reconstruction formula} holds for frames in Hilbert $C^*$-modules:
\begin{equation}\label{reconstruction}
x=SS^{-1}x=\sum_{j\in J}\langle S^{-1}x, x_j\rangle x_j=\sum_{j\in J}\langle x, S^{-1} x_j\rangle x_j,\quad x\in\h.
\end{equation}
We call $\{S^{-1}x_j\}_{j\in J}$ the {\it  canonical dual frame} of $\{x_j\}_{j\in J}$. Also if  $\{x_j\}_{j\in J}$ and $\{y_j\}_{j\in J}$ are frames and
$x=\sum_{j\in J}\langle x,y_j\rangle x_j,$ for all $x\in \h$, then  $\{x_j\}_{j\in J}$ and $\{y_j\}_{j\in J}$ are called {\it dual frames}. \\

Following \cite{frank}, 
a frame  $\{x_j\}_{j\in J}$ for a Hilbert ${\cal A}$-module $\h$ is  a {\it Riesz basis} if  the following two conditions are satisfied: 
\bei\item[(i)] $x_j\neq 0$ for all $j\in J$.
\item[(ii)] If $\sum_{j\in S}a_jx_j=0$  for a finite set of coefficients $\{a_j\}_{j\in S}\subset {\cal A}$, $S\subseteq  J$, then  $a_jx_j=0$  for every $j\in S$.
\eni
 It is proved in  \cite{han-jing} that a Riesz basis in a Hilbert $C^*$-module  may have many dual frames and it may even admit two different dual frames both of which are Riesz bases. This shows that the definition of Riesz bases in Hilbert $C^*$-modules is not the analog of Riesz basis in Hilbert spaces. Later, in \cite{khosravi} modular Riesz bases were introduced:

\bd
A sequence $\{x_j\}_{j\in J}$ in  Hilbert ${\cal A}$-module $\h$ is called a {\it modular Riesz basis} for $\h$, if there exists an invertible $\A$-linear and adjointable operator $U:l^2( J,{\cal A})\to \h$ such that $Ue_j=x_j$ for each $j\in J$, where $\{e_j\}_{j\in J}=\{(\delta_{ij}1_{\cal A})_{i\in J}\}_{j\in J}$ is the standard orthonormal basis of $l^2( J,{\cal A})$.
\ed

It is proved in \cite{ostad} that a sequence $\{x_i\}_{j\in J}\subset\h$ is a modular Riesz basis  if and only if
		$\{x_j\}_{j\in J}$ is a set of generators for $\h$ (as a Banach ${\cal A}$-module) and there exist $C,D>0$ such that for every finite sequence $\{a_i\}_{i\in S}$ in ${\cal A}$,
		\bee\label{0505a} C\left\|\sum_{i\in S} |a_i|^2\right\|\leq\left\|\sum_{i\in S} a_ix_i\right\|^2\leq D\left\|\sum_{i\in S} |a_i|^2\right\|.\ene
This gives a characterization of modular Riesz bases in terms of the analysis operator. As a consequence of  \eqref{0505a}, it is proved in \cite{ostad} that the set of all modular Riesz bases coincides with the set of all frames that have a unique dual frame which is a modular Riesz basis. Also it is straightforward from \eqref{0505a}, that the analysis operator of a modular Riesz basis is bijective and that an infinite series $\suj a_j x_j $ is convergent if and only if $\{a_j\}_{j\in J}\in\ell^2(J,\A)$.

\section{Modular Riesz bases versus Riesz bases}
It is well-known that in a Hilbert space $\h$, Riesz bases, $w$-independent frames, and exact frames are different names for the same class of sequences. We will now define the analogue version of the mentioned sequences  in Hilbert $C^*$-modules, and consider their interrelations.  Given a $C^*$-algebra $\A$ and  a Hilbert $\A$-module $\h$, let
\bes {\A}\text{-span}\{x_j\}_{j\in J}:=\left\{\sum_{j\in S}a_jx_j: a_j\in {\cal A}, \text { and }S\subset J\text{ is a finite subset }\right\}.\ens
\bd\label{0505b}
A sequence $\{x_j\}_{j\in J}\subset\h$ is said to be:
\bei\item[(i)]  $\omega$-independent, if whenever $\sum_{j\in J}a_jx_j$ is convergent and equal to zero for some  $\{a_j\}_{j\in J}\subset {\cal A}$, then necessarily $a_j=0$ for all $j\in J$.
		\item[(ii)]  biorthogonal with $\{y_j\}_{j\in J}\subset \h$, if $\langle x_i,y_j \rangle =\delta_{ij}1_{\cal A}$ for all $i,j\in J.$
		\item[(iii)]   exact frame, if $\{x_j\}_{j\in J}$ is a frame and for every $\ell\in J$,  $\{x_j\}_{j\neq \ell}$ cease to be a frame for $\h$. 
\eni\ed

The following result generalizes Theorem 7.7.1  in \cite{CB} to Hilbert $C^*$-modules and  shows  the connection between the sequences defined in part $(i)$ and $(ii)$ of Definition \ref{0505b} and modular Riesz bases. 
\bt\label{marzi}
	Let $\{x_j\}_{j\in J}$ be a frame in a Hilbert ${\cal A}$-module $\h$. Then the following statements are equivalent:
	\bei
		\item[(i)] $\{x_j\}_{j\in J}$ is a modular Riesz basis,
		\item [(ii)]$\{x_j\}_{j\in J}$ is $\omega$-independent,
		\item[(iii)] $\{x_j\}_{j\in J}$ and its canonical dual $\{S^{-1} x_j\}_{j\in J}$ are biorthogonal,
		\item[(iv)] $\{x_j\}_{j\in J}$ has a biorthogonal sequence,
	\eni
\et
\bp
$(i)\Rightarrow (ii) $ If $\{x_j\}_{j\in J}$ is a modular Riesz basis, then $T^*:l^2( J, {\cal A})\rightarrow \h$ is bijective by \eqref{0505a}. Assume that $\sum_{j\in J}a_j x_j=0$ for some $\{a_j\}_{j\in J}\subseteq {\cal A}$. Since $\sum_{j\in J}a_jx_j$ is convergent,  $\{a_j\}_{j\in J}\in l^2( J, {\cal A})$. Therefore $\{a_j\}_{j\in J}\in Ker T^*$. Since $T^*$ is injective, we have $\{a_j\}_{j\in J}=0$. \\
$(ii)\Rightarrow (iii) $ Consider the canonical dual frame $\{S^{-1}x_j\}_{j\in J}$, where $S$ is the frame operator of $\{x_j\}_{j\in J}$.
We have $x_\ell=\sum_{j\in J}\langle S^{-1}x_\ell,x_j\rangle x_j$ for all $\ell\in J$.
Hence
$$\sum_{j\neq \ell}\langle S^{-1}x_\ell,x_j\rangle x_j+(\langle S^{-1}x_\ell,x_\ell\rangle -1_{\cal A}) x_\ell=0.$$
Since $\{x_j\}_{j\in J}$ is $\omega$-independent, we have
$$\langle S^{-1}x_\ell,x_\ell\rangle =1_{\cal A},\qquad \langle S^{-1}x_\ell,x_j\rangle=0\quad for~each~j\neq \ell.$$
$(iii)\Rightarrow (iv)$ Clear.\\
$(iv)\Rightarrow (i)$ Let $T^*$ be the synthesis operator for $\{x_j\}_{j\in J}$ and let $\{a_j\}_{j\in J}\in Ker T^*$ and  $\{y_j\}_{j\in J}$ be biorthogonal to $\{x_j\}_{j\in J}$. Then for every $\ell\in J$, 
\bes a_\ell = \la a_\ell  x_\ell , y_\ell \ra = \la \suj a_j x_j , y_\ell \ra = 0. \ens 
Therefore the synthesis operator is bijective.
 \ep

\bc
	Let $\h$ be a Hilbert ${\cal A}$-module $\h$ over a unital $C^*$-algebra ${\cal A}$ and let $\{x_j\}_{j\in J}$ be a frame in $\h$ with the canonical dual frame $\{S^{-1}x_j\}_{j\in J}$. If $\langle x_j,S^{-1}x_j\rangle=1_{\cal A}$ for each $j\in J$, then $\{x_j\}_{j\in J}$ is a modular Riesz basis.
\ec
\bp  By the reconstruction formula, for each $j\in J$ we have
$$x_j=\sum_{\ell\in J}\langle x_j,S^{-1}x_\ell\rangle x_\ell=x_j+\sum_{\ell\neq j}\langle x_j,S^{-1}x_\ell\rangle x_\ell.$$
Thus $\sum_{\ell\neq j}\langle x_j,S^{-1}x_\ell\rangle x_\ell=0$. Applying $S^{-1}$, we obtain $\sum_{\ell\neq j}\langle x_j,S^{-1}x_\ell\rangle S^{-1}x_\ell=0$. From this relation, we conclude that
\bee\label{170819}\sum_{\ell\neq j}\langle x_j,S^{-1}x_\ell\rangle \langle S^{-1}x_\ell,x_j\rangle=0, \quad \text{ for all }j\in J.\ene
Since the element $\langle x_j,S^{-1}x_\ell\rangle \langle S^{-1}x_\ell,x_j\rangle = \langle x_j,S^{-1}x_\ell\rangle\langle x_j,S^{-1}x_\ell\rangle^*$ is  positive  for all $j,\ell\in J$, 
\eqref{170819} implies that 
\bes \| \langle x_j,S^{-1}x_\ell\rangle\|^2 = \| \langle x_j,S^{-1}x_\ell\rangle\langle x_j,S^{-1}x_\ell\rangle^*\| = 0\quad\text{ for all }\ell\in J\setminus \{j\}.\ens This means that $\{x_j\}_{j\in J}$ and its canonical dual frame are biorthogonal. Therefore by Theorem \ref{marzi}, $\{x_j\}_{j\in J}$ is a modular Riesz basis. \ep

In a Hilbert space, a frame is exact if and only if it is a Riesz basis, \cite[Theorem 5.5.4]{CB}.
We will now analyse the relationship between the exact frames and Riesz bases in  Hilbert $C^*$-modules; in particular the result with show that the situation is different compare to the Hilbert space case. We will need the following lemma, which  yields a characterization for exact frames in Hilbert $C^*$-modules, see \cite{jing} for the proof. 

\bl\label{jing}
	Let $\{x_j\}_{j\in J}$ be a frame for a Hilbert ${\cal A}$-module  $\h$ and let  $1_{\cal A}$ be the identity element of  ${\cal A}$. For each $\ell\in J$, the sequence $\{x_j\}_{j\neq \ell}$ is a frame for $\h$ if and only if $1_{\cal A}-\langle x_\ell,S^{-1}x_\ell\rangle$ is invertible in ${\cal A}$.
\el
The following result shows that not only modular Riesz basis but also the Riesz bases are indeed exact frames.
\begin{proposition}
	Let $\h$ be a Hilbert ${\cal A}$-module $\h$ over a unital $C^*$-algebra ${\cal A}$. 
If $\{x_j\}_{j\in J}$ is a Riesz basis then it is an exact frame.
\end{proposition}
\bp  Assume that $\{x_j\}_{j\in J}$ is a Riesz basis and that there exists some $\ell\in J$ such that $\{x_j\}_{j\neq \ell}$ is a frame for $\h$. Then there exists $\{a_j\}_{j\in J\setminus\{\ell\}}\subset {\cal A}$ such that $x_\ell=\sum_{j\neq \ell}a_jx_j$. Hence letting $a_\ell=-1_{\cal A}$, we have $\sum_{j\in  J}a_jx_j=0$. This leads to a contradiction with the definition of a Riesz basis. 
\ep


Unfortunately, as we will see in the next example,  the set of modular Riesz bases is smaller than the set of exact frames. 
\bex
	Consider the Banach space of all bounded sequences $\ell^\infty(\mn)$. This space forms a $C^*$-algebra with respect to operations  \bes \{u_j\}_{j\in\Bbb N}.\{v_j\}_{j\in\Bbb N}=\{u_jv_j\}_{j\in\Bbb N}, \quad \{u_j\}_{j\in\Bbb N}^*=\{\bar{u}_j\}_{j\in\Bbb N}.\ens
	  Let $\h=c_0$ denote the Banach space of all sequences vanishing at infinity. Then $\h$ is  a Hilbert ${\cal A}$-module with inner product 
	$$\langle \,\{u_j\}_{j\in\Bbb N}, \{v_j\}_{j\in\Bbb N}\rangle=\{u_j\bar{v_j}\}_{j\in\Bbb N},\quad \{u_j\}_{j\in\Bbb N},\{v_j\}_{j\in\Bbb N}\in\h.$$   Let $\delta_j\in \h$ be the sequence in $\h$ that takes the value 1 at the jth coordinate and 0 everywhere else. Note that the sequence $\{\delta_j\}_{j\in J}$ is a Parseval frame, since for every $u=\{u_j\}_{j\in J}\in \h$
	\bes \left\| \suj \la u,\delta_j\ra \la \delta_j , u \ra  \right\| &= &\left\| \suj |u_j|^2\delta_j  \right\| = \left\| \{|u_j|^2\}_{j\in J} \right\|\\
	&=&\sup |u_j|^2 = (\sup |u_j|)^2 = \| u\|^2.
	\ens
	Thus the frame operator is $S=Id_\h$ and  for each $j\in\Bbb N$,
	$$\langle \delta_j,S^{-1}\delta_j\rangle=\langle \delta_j, \delta_j\rangle= \delta_j.$$
	Therefore the sequence $1_{\cal A}-\langle \delta_j,S^{-1}\delta_j\rangle$ takes the value 0 at the jth coordinate and 1 everywhere else.  This shows  that $\langle \delta_j,S^{-1}\delta_j\rangle\neq 1_{\cal A}$ and that also  $1_{\cal A}-\langle \delta_j,S^{-1}\delta_j\rangle$ is not invertible in ${\cal A}$. By Lemma \ref{jing} we conclude that $\{\delta_j\}_{j\in\Bbb N}$  is an exact frame. On the other hand Theorem \ref{marzi} implies that $\{\delta_j\}_{j\in\Bbb N}$  is not a modular Riesz basis. Moreover, the set  $\{\delta_j\}_{j\in\Bbb N}$ is indeed a Riesz basis for $\h$. To see this, let $a_j=\{a_{ij}\}_{i\in\mn}\in \A$ and $S\subset \mn$ be a finite set and assume that $\sum_{j\in S} a_j \delta_j =0$.  We have 
	\bes \sum_{j\in S} a_j \delta_j =  \sum_{j\in S} \{ a_{ij}\delta_{ij}\}_{i\in\mn} =    \{ \sum_{j\in S} a_{ij}\delta_{ij}\}_{i\in\mn}=0. \ens
	This implies that $a_{jj}=0$ for all $j\in S$ and therefore $a_j\delta_j = 0$. This shows that $\{\delta_j\}_{j\in\Bbb N} $ is a Riesz basis for $\h$. \ep
\enx
\subsection*{Acknowledgments} 
After acceptance of the paper the author became award of the paper \cite{AB} which has some
overlap with the results. Indeed, part of the results in Theorem~\ref{marzi} also follows from
Cor. 4.20 in \cite{AB}. The author would like to thank Professor Arambasic for pointing this out.











\begin{thebibliography}{99}
	
	\bibitem{arambasic}
	L. Aramba\v{s}i\'{c}, {\it On frames for countably generated Hilbert $C^*$-modules}, Proc. Amer. Math. Soc. 135 (2007) 469--478.
	\bibitem{AB}
	L. Aramba\v{s}i\'{c}, and  D. Baki\'{c}. {\it Frames and outer frames for Hilbert-modules.} Linear and multilinear algebra 65.2 (2017) 381--431.
	
	\bibitem{bakic}
	D. Baki\'{c}, {\it Weak frames in Hilbert $C^*$-modules with application in Gabor analysis.} arXiv preprint arXiv:1903.01952 (2019).
	
	
\bibitem{CB} O. Christensen:
{\it An introduction to frames and Riesz bases,} Second expanded
edition. Birkh\"auser (2016).
	
	
	\bibitem{larson2}
	M. Frank and D. R. Larson, {\it A module frame concept for Hilbert $C^*$-modules}, Contemporary Math., 247 (1999) 207--234.
	
	\bibitem{frank}
	M. Frank and D. R. Larson, {\it Frames in Hilbert $C^*$-modules and $C^*$-algebras}, J. Operator Theory 48 (2002) 273--314.
%
	\bibitem{han-jing}
	D. Han, W. Jing, D. Larson and R. Mohapatra, {\it Riesz bases and their dual modular frames in Hilbert $C^*$-modules}, J. Math. Anal.
	Appl. 343 (2008) 246--256.

	
	\bibitem{jing}
	W. Jing, {\it Frames in Hilbert $C^*$-modules}, PhD thesis, University of Central Florida, 2006.
	
	
	\bibitem{kas}
	I. Kaplansky, {\it Algebras of type I}, Ann. Math. 56 (1952) 460--472.
	\bibitem{khosravi}
	A. Khosravi, B. Khosravi , {\it Fusion frames and g-frames in Hilbert $C^*$-modules}, Int. J. Wavelet, Multiresolution and Inf. Processing 6 (2008) 433--446.
	\bibitem{ostad}
	A. Khosravi, B. Khosravi, {\it g-frames and modular Riesz bases in Hilbert $C^*$-modules}, Int. J. Wavelets, Multiresolution and Inf. Processing. 10 (2012) 1250013.
	\bibitem{moosa}
	A. Khosravi, K. Musazadeh, {\it Fusion frames and g-frames}, J. Math. Anal. Appl. 342 (2008) 1068--1083.
	
	\bibitem{lance}
	E. C. Lance, {\it Hilbert $C^*$-modules: A Toolkit for Operator Algebraists}, London Math. Soc. Lecture Note Ser. vol. 210, Cambridge Univ. Press, 1995.
	\bibitem{luef}
	A. Austad, M. S. Jakobsen, and F. Luef, {\it Gabor Duality Theory for Morita Equivalent $ C^* $-algebras.} arXiv preprint arXiv:1905.01889 (2019).
	
	
	
	\bibitem{wegge}
	N. Wegge-Olsen, {\it K-theory and $C^*$-algebras - A Friendly Approach}, Oxford Univ. Press, Oxford, England, 1993.
\end{thebibliography}

\vspace{0.85cc} {\small \noindent Marzieh Hasannasab,
	Technical University of Kaiserslautern\\
	Paul-Ehrlich Stra\ss e Geb\"aude 31, 67663 Kaiserslautern,  Germany\\
	hasannas@mathematik.uni-kl.de\\\\ }

\end{document}